\documentclass{amsart}
\usepackage{amstext,amssymb,amsthm,amsopn,newlfont,graphpap,graphics,graphicx,mathrsfs,enumitem}
\allowdisplaybreaks
\usepackage[parfill]{parskip}
\usepackage[noadjust]{cite}
\usepackage{epigraph}
\usepackage[colorlinks=true,
            linkcolor=red,
            urlcolor=blue,
            citecolor=magenta]{hyperref}
\usepackage{color}
\usepackage{mathrsfs}
\allowdisplaybreaks
\theoremstyle{plain}

\newtheorem*{thm*}{Theorem}

\newtheorem*{prop*}{Proposition}

\newtheorem*{cor*}{Corollary}

\newtheorem*{lem*}{Lemma}
\theoremstyle{definition}

\newtheorem*{defn*}{Definition}

\newtheorem*{exmps*}{Examples}

\newtheorem*{exmp*}{Example}

\newtheorem*{exerc*}{Exercise}

\newtheorem{rems}{Remarks}[section]
\newtheorem*{rems*}{Remarks}

\newtheorem*{rem*}{Remark}
\newcommand{\N}{{\mathbb N}}
\newcommand{\C}{{\mathbb C}}

\newcommand{\emps}{\emptyset}

\renewcommand{\bar}{\overline}

\numberwithin{equation}{section}

\begin{document}
\title[On unbounded linear operators with an arbitrary spectrum]
{An observation on unbounded linear operators with an arbitrary spectrum}
\author[Marat V. Markin]{Marat V. Markin}
\address{
Department of Mathematics\newline
California State University, Fresno\newline
5245 N. Backer Avenue, M/S PB 108\newline
Fresno, CA 93740-8001, USA
}
\email{mmarkin@csufresno.edu}
\subjclass{Primary 47A10; Secondary 47A08, 47B15}
\keywords{Bounded linear operator, closed  linear operator, normal operator, spectrum}
\begin{abstract}
We furnish a simple way of constructing an unbounded closed linear operator in a complex Banach space, whose spectrum is an arbitrary nonempty closed, in particular compact, subset of the complex plane.
\end{abstract}
\maketitle

\section[Introduction]{Introduction}

By \textit{Gelfand's spectral radius theorem} (see, e.g., \cite{Markin2020EOT}), the spectrum of a bounded linear operator on a complex Banach space is a nonempty compact subset of the complex plane. Conversely, an arbitrary nonempty compact set in $\C$ is the spectrum of a such an operator (cf. \cite[Theorem $5.7$]{Markin2020EOT}).

While the spectrum of an unbounded closed linear operator in a complex Banach space is known to be a closed subset of the complex plane, including $\emps$ and $\C$, (see, e.g., \cite{Dun-SchI,Markin2020EOT}), the natural question is whether an arbitrary closed set in $\C$ is the spectrum of a such an operator. Of particular interest are unbounded closed linear operators with nonempty compact spectra.

We give a simple way of constructing an unbounded closed linear operator in a complex Banach space, whose spectrum is any given nonempty closed, in particular compact, subset of the complex plane.

\section[Preliminaries]{Preliminaries}

For a closed linear operator $A$ in a complex Banach space $X$, the set
\[
\rho(A):=\left\{ \lambda\in \C \,\middle|\, \exists\, R(\lambda,A):=(A-\lambda I)^{-1}\in L(X) \right\}
\]
($I$ is the \textit{identity operator} on $X$, $L(X)$ is the space of bounded linear operators on $X$) and its complement $\sigma(A):=\C\setminus \rho(A)$ are called the operator's \textit{resolvent set} and \textit{spectrum}, respectively.

The spectrum is partitioned into three pairwise disjoint subsets, $\sigma_p(A)$, $\sigma_c(A)$, and $\sigma_r(A)$, called the \textit{point}, \textit{continuous}, and \textit{residual spectrum} of $A$, respectively, as follows:
\begin{equation*}
\begin{split}
& \sigma_p(A):=\left\{\lambda\in \C \,\middle|\,A-\lambda I\ \text{is \textit{not injective}, i.e., $\lambda$ is an \textit{eigenvalue} of $A$} \right\},\\
& \sigma_c(A):=\left\{\lambda\in \C \,\middle|\,A-\lambda I\ \text{is \textit{injective},
\textit{not surjective}, and $\overline{R(A-\lambda I)}=X$} \right\},\\
& \sigma_r(A):=\left\{\lambda\in \C \,\middle|\,A-\lambda I\ \text{is \textit{injective} and $\overline{R(A-\lambda I)}\neq X$} \right\}
\end{split}
\end{equation*}
($R(\cdot)$ is the \textit{range} of an operator and $\overline{\cdot}$ is the \textit{closure} of a set) (see, e.g., \cite{Dun-SchI,Markin2020EOT}).

In the complex Hilbert space $l_2$ of square-summable sequences, the operator of multiplication by a sequence
$\left(m_n\right)_{n\in \N}\in \C^\N$ 
\[
M\left(x_n\right)_{n\in \N}:=\left(m_nx_n\right)_{n\in \N}
\]
with maximal domain
\[
D(M):=\left\{ \left(x_n\right)_{n\in \N}\in l_2\,\middle|\, \left(m_nx_n\right)_{n\in \N}\in l_2 \right\}
\] 
is a (densely defined) closed linear operator, which is bounded \textit{iff} the multiplier sequence is bounded, i.e.,
\[
\sup_{n\in \N}|m_n|<\infty.
\] 

In fact, the operator $M$ is \textit{normal} (see, e.g., \cite{Dun-SchII}) and
\[
\sigma(M)=\bar{\left\{m_n\right\}}
\]
with
\[
\sigma_p(M)=\left\{m_n\right\},\ 
\sigma_c(M)=\bar{\left\{m_n\right\}}\setminus \left\{m_n\right\},
\ \text{and}\
\sigma_r(M)=\emps,
\]
where $\left\{m_n\right\}$ is the set of values of $(m_n)_{n\in\N}$
\cite{Markin2020EOT}.

Thus, by choosing such multiplier sequence $(m_n)_{n\in\N}$ that $\{m_n\}_{n\in\N}$ is a countable dense subset of an arbitrary nonempty closed set $\sigma$ in the complex plane, we obtain the multiplication operator $M$ in $l_2$ with
\[
\sigma(M)=\sigma,
\]
the normal operator $M$ being bounded whenever the set $\sigma$ is compact (see, e.g., \cite{Dun-SchII}).

In the complex Banach space $Y:=L_p(0,1)$ ($1\le p<\infty$) or $Y:=C[0,1]$, the latter equipped with the maximum norm
\[
C[0,1]\ni x\mapsto \|x\|_\infty:=\max_{0\le t\le 1}|x(t)|,
\]
the differentiation operator
\[
Dx:=x'
\]
with domain 
\begin{equation*}
D(D):=\left\{ x\in  L_p(0,1)\,\middle|\, x(\cdot)\in AC[a,b],\
x'\in L_p(0,1),\ x(0)=0\right\}
\end{equation*}
or
\[
D(D):=\left\{ x\in  C^1[0,1]\,\middle|\, x(0)=0\right\},
\]
respectively, is an unbounded closed linear operator with
\[
\sigma(D)=\emps
\]
and
\[
[R(\lambda,D)y](t)=\int_0^t e^{\lambda(t-s)}y(s)\,ds,\ \lambda\in\C,y\in Y,
\]
(cf. \cite{Markin2020EOT}), the operator being densely defined in
$L_p(0,1)$ ($1\le p<\infty$) but not in $C[0,1]$.

\section{Unbounded Linear Operator with an Arbitrary Spectrum}

Let $\sigma$ be arbitrary nonempty closed set in the complex plane,
$M$ be the multiplication operator in the complex space $l_2$ relative to a multiplier sequence $\left(m_n\right)_{n\in \N}\in \C^\N$ such that $\left\{m_n\right\}$ is a countable dense subset of $\sigma$, and hence, 
\[
\sigma(M)=\bar{\left\{m_n\right\}}=\sigma,
\]
and $D$ be the differentiation operator in the complex space $L_p(0,1)$ ($1\le p<\infty$) or $C[0,1]$ with $\sigma(D)=\emps$ (see Preliminaries).

In the complex Banach space $X\oplus Y$, where $X:=l_2$
and $Y:=L_p(0,1)$ ($1\le p<\infty$) or $Y:=C[0,1]$,
equipped with the norm
\[
X\oplus Y\in (x,y)\mapsto \|(x,y)\|:=\|x\|_X+\|y\|_Y,
\]
the operator matrix
\[
\begin{bmatrix}
M&0\\
0&D\\
\end{bmatrix}
\]
defines an unbounded closed linear operator $A$, whose domain 
\[
D(A):=D(M)\oplus D(D)
\]
is dense in $X\oplus Y$ whenever $D(D)$ is dense in $Y$, with
\[
\sigma(A)=\sigma(M)\cup \sigma(D)=\sigma(M)=\bar{\left\{m_n\right\}}=\sigma
\]
and, as is easily seen,
\[
\sigma_p(A)=\sigma_p(M)=\left\{m_n\right\},\
\sigma_c(A)=\sigma_c(M)=\bar{\left\{m_n\right\}}\setminus \left\{m_n\right\},\ \text{and}\
\sigma_r(A)=\emps
\]
(see Preliminaries).

\begin{rems}\
\begin{itemize}
\item Thus, the spectrum of an unbounded closed linear operator in a complex Banach space can be an arbitrary closed subset of the complex plane, including $\emps$ and $\C$, (cf. \cite[Theorem $5.7$]{Markin2020EOT}). 
\item If, in the prior example, $Y:=L_2(0,1)$, 
$X\oplus Y$ is a Hilbert space relative to the inner product
\[
X\oplus Y\in (x_1,y_1),(x_2,y_2)\mapsto \langle(x_1,y_1),(x_2,y_2)\rangle:=\langle x_1,x_2\rangle_X
+\langle y_1,y_2\rangle_Y.
\]
\end{itemize}
\end{rems}


\end{document}